\documentclass[a4paper,12pt]{amsart}

\usepackage{amsmath,amssymb}

\textwidth=\paperwidth \advance\textwidth by-2in
\advance\oddsidemargin by-1in
\evensidemargin=\oddsidemargin

\def\Z{\mathbb{Z}}
\def\Q{\mathbb{Q}}
\def\R{\mathbb{R}}

\def\Hom{\operatorname{Hom}}

\theoremstyle{plain}
\newtheorem{thm}{Theorem}
\newtheorem{prop}[thm]{Proposition}

\newtheorem{lem}[thm]{Lemma}
\theoremstyle{definition}

\newtheorem*{remark}{Remark}

\begin{document}

\title[Murasugi polynomial] {A characterization of the Murasugi
  polynomial of an equivariant slice knot}

\author{Jae Choon Cha}

\address{Information and Communications University, Daejeon 305--714,
  Korea}

\email{jccha@icu.ac.kr}

\def\subjclassname{\textup{2000} Mathematics Subject Classification}
\expandafter\let\csname subjclassname@1991\endcsname=\subjclassname
\expandafter\let\csname subjclassname@2000\endcsname=\subjclassname


\begin{abstract}
  We characterize the Murasugi polynomial of an equivariant slice knot
  by proving a conjecture of J. Davis and S. Naik.
\end{abstract}

\maketitle

A knot $K$ in $S^3$ is called \emph{periodic of period $p$} if there
is an orientation preserving action of $\Z/p$ on $S^3$ which preserves
$K$ setwise and its fixed point set $A$ is a circle disjoint to~$K$.
$A$~is called the \emph{axis}.  Two periodic knots $K_0$ and $K_1$ of
period $p$ are called \emph{equivariantly concordant} if there is an
action of $\Z/p$ on $S^3\times[0,1]$ such that $K_i\times i$ is
periodic via its restriction on $S^3\times i$, $i=0,1$, and there is a
locally flat submanifold in $S^3\times[0,1]$ which is preserved by the
action, homeomorphic to $S^1\times [0,1]$, and bounded by $(K_0 \times
0) \cup -(K_1\times 1)$.  The unknotted circle $S^1\times 0 \subset
\partial (D^2\times D^2)=S^3$ can be viewed as a periodic knot via the
$(2\pi/p)$-rotation on the first $D^2$ factor.  If a periodic knot $K$
is equivariantly concordant to it, $K$ is called an \emph{equivariant
  slice knot}.

There are several known obstructions for a periodic knot $K$ to being
an equivariant slice knot.  Some of them are obtained from invariants
of~$K$.  In~\cite{Naik:1997-1}, Naik used the Alexander polynomial and
metabolizers of the Seifert form of~$K$.  She also showed that certain
Casson-Gordon invariants of $K$ must vanish if $K$ is equivariant
slice.  In~\cite{Choi-Ko-Song:1998-1}, Choi, Ko, and Song defined an
obstruction from a Seifert matrix of~$K$.

Further obstructions are obtained by considering the quotient link.
Given a periodic knot $K$ with axis $A$, the orbit space of the
$(\Z/p)$-action is again $S^3$ by the Smith conjecture, and the images
$\bar A$ and $\bar K$ of $A$ and $K$ under the quotient map form a
two-component link which is called the \emph{quotient link}.  It
contains all the essential information on the periodic knot.
In~\cite{Cha-Ko:1999-2}, Ko and the author developed an obstruction
for $K$ to being an equivariant slice knot from knots obtained by
surgery on the quotient link.  In particular, their Casson-Gordon
torsion invariant was used to construct an example of a
non-equivariant-slice knot which cannot be detected by other
invariants.

Recently, in~\cite{Davis-Naik:2002-1}, Davis and Naik have studied the
Murasugi polynomial $\Delta_{\Z/p}(g,t)$ of a periodic knot~$K$, which
is the image of the Alexander polynomial of the quotient link under
the projection $\Z[\Z\times \Z] \to \Z[\Z/p\times \Z]$.  Here $g$ and
$t$ are generators of $\Z/p$ and $\Z$ corresponding to the components
$\bar A$ and $\bar K$, respectively.  They proved the following
realization theorem of the Murasugi polynomial of an equivariant
slice knot:

\begin{thm}[Davis-Naik]
  For any $a(g,t)\in \Z[\Z/p\times \Z]$ such that $a(g,1)=1$, there is
  an equivariant slice knot $K$ with Murasugi polynomial
  $\Delta_{\Z/p}(g,t)=a(g,t)a(g^{-1},t^{-1})$.
\end{thm}

In fact, their knot $K$ is an equivariant \emph{ribbon} knot, which is
a specialization of an equivariant slice knot.

They conjectured that the converse is true, i.e., the Murasugi
polynomial $\Delta_{\Z/p}(g,t)$ of every equivariant slice knot is of
the above form.  Some related results have been revealed.
In~\cite{Davis-Naik:2002-1}, by interpreting the Murasugi polynomial
as the Reidemeister torsion, Davis and Naik proved that if $K$ an
equivariant slice knot, then
$$
\Delta_{\Z/p}(g,t)b(g,t)b(g^{-1},t^{-1})=a(g,t)a(g^{-1},t^{-1})
$$
up to $\pm g^r t^s$, for some $a, b \in \Z[\Z/p\times \Z]$ such
that $a(g,1)=1=b(g,1)$.  So the question becomes whether $b=1$ in this
result.  They also showed that $b=1$ for an equivariant ribbon knot.
In~\cite{Hillman:2002-1}, Hillman proved that if $K$ is an equivariant
slice knot, then $\Delta_{\Z/p}(g,t)=a(g,t)a(g^{-1},t^{-1})$, up to
units, over $\Q[\Z/p\times \Z]$.

The goal of this memo is to prove the Davis-Naik conjecture:
\begin{thm}\label{thm:main}
  If $K$ is an equivariant slice knot, then
  $\Delta_{\Z/p}(g,t)=a(g,t)a(g^{-1},t^{-1})$, up to $\pm g^r t^s$,
  for some $a(g,t) \in \Z[\Z/p\times \Z]$ such that $a(g,1)=1$.
\end{thm}

Combined with the above result of Davis and Naik, it characterizes the
Murasugi polynomial of an equivariant slice knot.

%

For the proof of Theorem~\ref{thm:main} we use the Blanchfield forms
of links which have a well developed theory in the literature.  Our
arguments are based on the ideas and results of
Blanchfield~\cite{Blanchfield:1957-1}, Hillman~\cite{Hillman:2002-1},
Levine~\cite{Levine:1982-1}.  We will focus on a special case where we
have sharpened versions of well known general results.
Theorem~\ref{thm:main} will follow from results of this special case.

Let $\Lambda$ be the group ring $\Z[\Z\times \Z]$ which is identified
with the ring of Laurent polynomials in variables $x, y$, and let $Q$
be its quotient field.  For a two-component link $L$, we denote its
exterior by~$E_L$.  The abelianization map $\pi_1(E_L) \to \Z\times\Z$
sending meridians to the standard basis gives rise to a
$\Lambda$-coefficient system on~$E_L$.  For a $\Lambda$-module $M$,
denote its torsion submodule by~$tM$.  A~$\Lambda$-module is called
\emph{pseudozero} if its localization away from $\pi$ is zero for all
prime $\pi$ in~$\Lambda$.  (Since $\Lambda$ is a UFD, it agrees with
the standard general definition requiring that $\pi$ is of height
one.)  We denote by $\hat tM$ the quotient of $tM$ by its maximal
pseudozero submodule.  Then there is a non-degenerated sesquilinear
pairing
$$
\hat t H_1(E_L;\Lambda) \times \hat t H_1(E_L,\partial E_L;\Lambda)
\to Q/\Lambda
$$
due to Blanchfield~\cite{Blanchfield:1957-1}.  (We say $A\times B
\to C$ is non-degenerated if the adjoint maps $A\to \Hom(B,C)$ and $B
\to \Hom(A,C)$ are injective.)  It induces
$$
B_L\colon \hat t H_1(E_L;\Lambda) \times \hat t H_1(E_L;\Lambda)
\to Q/\Lambda
$$
which is not necessarily non-degenerated.  In general, in order to
obtain a non-degenerated pairing whose Witt class is a concordance
invariant of an arbitrary link, $B_L$ is localized by inverting an
appropriate multiplicative subset of $\Lambda$ (for example, see
Hillman~\cite{Hillman:2002-1}).

From now on we assume that a link $L$ has two components with
nontrivial linking number.  In this special case, we have the crucial
advantage that the \emph{unlocalized} Blanchfield form $B_L$ is
invariant under concordance.

\begin{lem}\label{lem:prop-of-link-with-nonzero-lk}
\mbox{}
\begin{enumerate}
\item $H_i(\partial E_L;\Lambda)=0$ for $i>0$.
\item $H_1(E_L;\Lambda)$ is a torsion $\Lambda$-module.
\end{enumerate}
\end{lem}

\begin{proof}
  (1) follows from that the $(\Z\times\Z)$-cover of $\partial E_L$
  consists of copies of $\R^2$.  (2) is a result of
  Levine~\cite[Theorem A, page 378]{Levine:1982-1}.  See also
  Hillman~\cite{Hillman:2002-1}.
\end{proof}

From this $H_1(E_L;\Lambda) = tH_1(E_L;\Lambda)= tH_1(E_L,\partial
E_L;\Lambda)$ and the above $B_L$ is non-degenerated.  Furthermore, a
standard argument shows the following result on the unlocalized
Blanchfield form:

\begin{lem}\label{lem:invariance-of-B_L}
  If $L$ and $L'$ are concordant, then $B=B_L \oplus (-B_{L'})$ is
  metabolic, i.e., there is a submodule $P$ in $\hat t
  H_1(E_L;\Lambda)\oplus \hat t H_1(E_{L'};\Lambda)$ such that
  $P^\perp =P$ with respect to~$B$.
\end{lem}

The argument of Hillman~\cite[page 37--38]{Hillman:2002-1} proves this.  In
fact, in~\cite{Hillman:2002-1}, he proved an analogue for the Blanchfield
form over a certain localized coefficient system~$S^{-1}\Lambda$.  The
advantage of his localization is that the following fact holds for
\emph{any} link: if $W$ is the exterior of a concordance between $L$ and
$L'$, then $H_1(\partial W;S^{-1}\Lambda) \cong
H_1(E_L;S^{-1}\Lambda)\oplus H_1(E_{L'};S^{-1}\Lambda)$.  In our case, it
also holds for the (unlocalized) $\Lambda$-coefficient homology modules,
since $H_1(\partial E_L;\Lambda)=H_1(\partial E_{L'};\Lambda)=0$ by
Lemma~\ref{lem:prop-of-link-with-nonzero-lk}.  From this it can be seen
that the argument in~\cite{Hillman:2002-1} works for the unlocalized
Blanchfield pairing~$B_L$.  We omit the details.

For a finitely generated $\Lambda$-module $M$, $\Delta(M)$ is defined
to be the greatest common divisor of $n\times n$ minors of a
presentation matrix of $M$ where $n$ is the number of generators of
the presentation.  It is known that $\Delta(M)$ of the underlying
module $M$ of a metabolic non-degenerated pairing is of the form
$f(x,y)f(x^{-1},y^{-1})$ (e.g.,
see~\cite{Blanchfield:1957-1,Hillman:2002-1}).  Thus the above lemma
implies
$$
\Delta(\hat t H_1(E_L;\Lambda)) \Delta(\hat t H_1(E_{L'};\Lambda)) =
f(x,y)f(x^{-1},y^{-1})
$$
up to units for some $f\in \Z[\Z\times \Z]$.  Recall that the
Alexander polynomial of $L$ is defined by
$\Delta_L(x,y)=\Delta(H_1(E_L;\Lambda))$.  Combining
Lemma~\ref{lem:prop-of-link-with-nonzero-lk} with a result of
Blanchfield~\cite[Theorem~4.7]{Blanchfield:1957-1} that
$\Delta(tM)=\Delta(\hat tM)$, it follows that
$$
\Delta(H_1(E_L;\Lambda)) = \Delta(tH_1(E_L;\Lambda)) = \Delta(\hat
tH_1(E_L;\Lambda)).
$$
From this we have

\begin{lem}\label{lem:invariance-of-Delta_L}
  If $L$ and $L'$ are concordant, $\Delta_L(x,y) \Delta_{L'}(x,y) =
  f(x,y)f(x^{-1},y^{-1})$, up to $\pm x^ry^s$, for some $f \in
  \Z[\Z\times\Z]$.
\end{lem}

If $L$ is a quotient link of an equivariant slice knot, then $L$ is
concordant to the Hopf link which has linking number one.  (Quotient
links of equivariantly concordant periodic knots are concordant in the
topologically locally flat category.)  Since the Hopf link has trivial
Blanchfield form and Alexnader polynomial, we obtain the following
consequence:

\begin{prop}
  If $L$ is a quotient link of an equivariant slice knot, then
  \begin{enumerate}
  \item The Blanchfield form $B_L$ is metabolic.
  \item The Alexander polynomial $\Delta_L(x,y)$ is of the form
    $f(x,y)f(x^{-1},y^{-1})$ for some $f \in \Z[\Z\times\Z]$
  \end{enumerate}
\end{prop}

In particular, reducing to $\Z[\Z/p\times \Z]$, we have
$\Delta_{\Z/p}(g,t)=a(g,t)a(g^{-1},t^{-1})$ over $\Z[\Z/p\times\Z]$ up
to $\pm g^r t^s$.  Since $\Delta_{\Z/p}(g,1)=1$, we may assume
$a(g,1)=1$.  This proves Theorem~\ref{thm:main}.

\begin{remark}
  In the above discussion we may avoid the use of $\hat t
  H_1(E_L;\Lambda)$ since Levine's result in ~\cite{Levine:1982-1}
  implies that $tH_1(E_L;\Lambda)=\hat t H_1(E_L;\Lambda)$ in our
  case.  But we still need to consider $\hat t H_1(E_L;\Lambda)$ to
  prove Lemma~\ref{lem:invariance-of-B_L} using the standard argument
  as in~\cite{Hillman:2002-1}.
\end{remark}

\begin{remark}
  Subsequent to this memo, Davis and Naik have found a different proof
  of Theorem~\ref{thm:main} using the Reidemeister torsion.
\end{remark}

\sloppy \bibliographystyle{amsplainabbrv}
\bibliography{research}

\end{document}